\documentclass[titlepage,twoside,12pt]{article}
\usepackage{amssymb}
\usepackage{amsfonts}
\begin{document}

{\bf \Large Self-Referential Definition of \\ \\ Orthogonality} \\ \\

{\bf Elem\'{e}r E Rosinger,~ Gusti van Zyl} \\ \\
Department of Mathematics \\
and Applied Mathematics \\
University of Pretoria \\
Pretoria \\
0002 South Africa \\
eerosinger@hotmail.com, gusti.vanzyl@up.ac.za \\ \\

{\bf Abstract} \\

There has for longer been an interest in finding equivalent conditions which define inner
product spaces, and the respective literature is considerable, see for instance Amir, which
lists 350 such results. \\
Here, in this tradition, an alternative definition of {\it orthogonality} is presented which
does {\it not} make use of any inner product. This definition, in the spirit of the recently
developed non-wellfounded set theory, is self-referential, or circulatory. \\ \\

{\bf 0. Preliminaries} \\

Let us note that a group $( G, \star)$ is defined by properties of its subgroups generated by no more than 3 elements.
Indeed, the axioms of a group are \\

(0.1) $~~~~~~ \forall~~ x, y \in G ~:~ x \star y \in G $ \\

(0.2) $~~~~~~ \forall~~ x, y, z \in G ~:~ x \star ( y \star z ) = ( x \star y ) \star z $ \\

(0.3) $~~~~~~ \exists~~ e \in G ~:~ \forall~~ x \in G ~:~ x \star e = e \star x = x $ \\

(0.4) $~~~~~~ \forall~~ x \in G ~:~ \exists~~ x\,' \in G ~:~ x \star x\,' = x\,' \star x = e $ \\

Clearly, vector spaces are also defined by properties of their vector subspaces generated by no more than 3 elements. \\

On the other hand, the conditions that a vector space $X$ be a normed space $( X, ||~|| )$ will only involve its 2
dimensional vector subspaces, namely \\

(0.5) $~~~~~~ \forall~~ x \in X ~:~ || x || \geq 0 $ \\

(0.6) $~~~~~~ \forall~~ x \in X ~:~ || x || = 0 ~~\Longleftrightarrow~~ x = 0 $ \\

(0.7) $~~~~~~ \forall~~ x, y \in X ~:~ || x + y || \leq || x || + || y || $ \\

Obviously, the same goes for finite dimensional Hilbert spaces, where the scalar product, and in particular,
orthogonality only involve 2-dimensional vector subspaces. \\

In the case of infinite dimensional vector spaces, just like with infinite dimensional Banach spaces, completeness is also
required, a condition which, of course, is no longer definable in terms of any finite dimensional vector subspaces. \\

It can in general be recalled the important role of the structure of finite dimensional vector spaces in the study of arbitrary normed spaces, as
illustrated by the earlier and more particular Khinchin, and the later Grothendieck inequalities. \\

In the sequel, we shall consider the issue of orthogonality, and do so without the need to deal with completeness. \\ \\

{\bf 1. Linear Independence} \\

Let $V$ be a vector space on $\mathbb{R}$, and $m \geq 2$, then we denote \\

(1.1) $~~~~~~ {\cal I}nd\,^m\,_V ~=~  \left \{~ ( a_1, ~.~.~.~ , a_m )~~
              \begin{array}{|l}
                     ~~ a_1, ~.~.~.~ , a_m \in V ~~\mbox{are linearly} \\
                     ~~~ \mbox {independent in}~~ V
              \end{array} ~\right \}$ \\

which is nonvoid, if and only if dim V $\geq$ m. \\

Further, for $( a_1, ~.~.~.~ , a_m ) \in {\cal I}nd\,^m\,_V$ we denote by \\

(1.2) $~~~~~~ span \,( a_1, ~.~.~.~ , a_m ) $ \\

the vector subspace in $V$ generated by $a_1, ~.~.~.~ , a_m$. Then for every m-dimensional vector subspace $E$ of $V$, we obviously we have \\

(1.3) $~~~~~~ \begin{array}{l}
                   \forall~~~ x \in E ~: \\ \\
                   \forall~~~ ( a_1, ~.~.~.~ , a_m ) \in {\cal I}nd\,^m\,_E ~: \\ \\
                   \exists~ ! ~~ \lambda_1 ( a_1, ~.~.~.~ , a_m, x ), ~.~.~.~ ,
                                \lambda_m ( a_1, ~.~.~.~ , a_m, x ) \in \mathbb{R} ~: \\ \\
                   ~~~~~~~ x ~=~ \Sigma_{1 \leq i \leq m}~ \lambda_i ( a_1, ~.~.~.~ , a_m, x )
                                             \, a_i
               \end{array} $ \\ \\

{\bf 2. Usual orthogonality based on inner product} \\

Here for convenience and for setting the notation, we review the usual case of orthogonality,
namely, that which is defined based on an {\it inner product}. Given any inner product
$< ~,~ >~ : V \times V \longrightarrow \mathbb{R}$, then for $m \geq 2$, we denote \\

(2.1) $~~~~~~ \begin{array}{l}
           {\cal O}rt\,^m\,_V ( < ~,~ > ) ~=~ \\ \\
           ~=~ \left \{~ ( a_1,.~.~.~. , a_m ) \in {\cal I}nd\,^m\,_V ~~
              \begin{array}{|l}~~
                         a_1,.~.~.~. , a_m ~~~\mbox{are orthogonal} \\
                         ~~~\mbox{with respect to}~ < ~,~ >
                    \end{array} ~\right \}
               \end{array} $ \\ \\

Now, if $( a_1, .~.~.~. , a_m ) \in {\cal O}rt\,^m\,_V ( < ~,~ > )$, then for any vector $x
\in span \, ( a_1,.~.~.~. , a_m )$, we have in a unique manner, see (1.3) \\

(2.2) $~~~~~~ x ~=~ \Sigma_{1 \leq i \leq m}~ \lambda_i ( a_1, ~.~.~.~ , a_m, x )\, a_i $ \\

where this time we have in addition the simple explicit formulae \\

(2.3) $~~~~~~ \lambda_i ( a_1, ~.~.~.~ , a_m, x ) ~=~ < a_i, x > / < a_i, a_i >,~~~
                                                                    1 \leq i \leq m $ \\

{\bf Remark 1.} \\

An important point to note is that each $\lambda_i( a_1, ~.~.~.~ , a_m, x )$ in (2.2), (2.3)
can in fact only depend on $x$ and $a_i$, but {\it not} on $a_j$, with $1 \leq j \leq m,~ j \neq
i$. Indeed, in view of (1.3), $x$ obviously does {\it not} depend on $a_1, ~.~.~.~ , a_m$.

\hfill $\Box$ \\

For the sake of clarity, we recall here that by inner product on $V$ we mean a function
$ < ~,~ > \,: V \times V \longrightarrow \mathbb{R}$, such that \\

(2.4)~~~ $  < ~,~ > ~~\mbox{is symmetric} $ \\

(2.5)~~~ $ < x, . > ~~\mbox{is linear on}~~ V $ \\

(2.6)~~~ $ < x, x > ~\geq~ 0 $ \\

(2.7)~~~ $ < x, x > ~=~ 0 ~~\Longleftrightarrow~~ x ~=~ 0 $ \\ \\

{\bf 3. Orthogonality without inner product} \\

In view of the above it may appear that one could give the following definition of {\it orthogonality without inner
product}, namely \\

{\bf Definition 1.} \\

If $m \geq 2$ and $( a_1, ~.~.~.~ , a_m ) \in {\cal I}nd\,^m\,_V$, then by definition \\ \\

(3.1) $\begin{array}{l}
       ~~~ a_1, ~.~.~.~ , a_m ~~~\mbox{are orthogonal} ~~~\Longleftrightarrow~~~ \\ \\ \\

       ~~~~~~~\Longleftrightarrow~~~
                  \left(~ \begin{array}{l}
                       \forall~~ x \in span \,( a_1, ~.~.~.~ , a_m ) ~: \\ \\
                       \forall~~ 1 \leq i \leq m ~: \\ \\
                       ~~~ \lambda_i ( a_1, ~.~.~.~ , a_m, x ) ~~\mbox{in}~~ (1.3)
                                ~~\mbox{does not} \\
                                ~~~\mbox{depend on}~ a_j\,,
                                 ~~\mbox{with}~~ 1 \leq j \leq m,~ j \neq i
                    \end{array} ~~\right )
        \end{array} $ \\ \\

Consequently, we denote \\

(3.2) $~~~ {\cal O}rt\,^m\,_V ~=~ \left \{~ ( a_1, ~.~.~.~ , a_m ) \in {\cal I}nd\,^m\,_V ~~
               \begin{array}{|l}~~
               a_1, ~.~.~.~ , a_m ~~~\mbox{are orthogonal} \\
               ~~~\mbox{in the sense of}~~ (3.1)
               \end{array} ~\right \}$ \\

{\bf Remark 2.} \\

One should note above that the meaning of "does not depend on" in (3.1) needs further clarification. Indeed, when we are
given an explicit expression, such as for instance in (2.3), then it is quite clear what "does not depend on" means, since the entities on which the
respective expressions are supposed not to depend simply do not appear in the respective expressions. \\

However, in (3.1) we are actually dealing with (1.3). And then, the meaning of \\

"~$\lambda_i ( a_1, ~.~.~.~ , a_m, x ) ~~\mbox{in}~~ (1.3) ~~\mbox{does not depend on}~ a_j\,, \\
                                                     ~~~~\mbox{with}~~ 1 \leq j \leq m,~ j \neq i$~" \\

is not so immediately obvious. \\ \\

{\bf 4. The two dimensional case} \\

Orthogonality, in its simplest nontrivial case, involves only two nonzero vectors, hence it is
a property which can already be formulated in 2-dimensional vector spaces. Let therefore $E$
be any 2-dimensional vector space on $\mathbb{R}$. Then the property in (1.3) becomes \\

(4.1) $~~~~~~\begin{array}{l}
     \exists~ !~~~ \lambda, \mu : {\cal I}nd\,^2\,_E \times E
                                ~~\longrightarrow~~ \mathbb{R} ~: \\ \\
     \forall ~~~~ ( a, b ) \in {\cal I}nd\,^2\,_E,~ x \in E ~: \\ \\
     ~~~~~~~ x ~=~ \lambda ( a, b, x )\, a + \mu ( a, b, x )\, b
           \end{array}$ \\ \\

Therefore, in view of (3.1) in the above Definition 1 which does {\it not} use inner product
for defining orthogonality, we have \\

{\bf Proposition 1.} \\

Given $( a, b ) \in {\cal I}nd\,^2\,_E$, then \\ \\

(4.2) $~~~ \begin{array}{l}
              \left ( \begin{array}{l}
                      a, b ~~~\mbox{are orthogonal} \\
                      \mbox{in the sense of}~~ (3.1)
                   \end{array} \right ) ~\Longleftrightarrow~ \\ \\ \\
              ~~~~~~~~~~~~~~~~\Longleftrightarrow~
              \left( \begin{array}{l}
                       \forall ~ x \in E ~: \\ \\
                       ~~(4.2.1)~~ \lambda ( a , b, x ) ~~\mbox{in}~~ (4.1) \\
                       ~~~~~~~~~~~~~\mbox{does not depend on}~ b \\ \\
                       ~~(4.2.2)~~ \mu ( a, b ,x ) ~~\mbox{in}~~ (4.1) \\
                       ~~~~~~~~~~~~~\mbox{does not depend on}~ a
                     \end{array} \right )
            \end{array} $ \\ \\

\hfill $\Box$ \\

Clearly, if we want to consider the concept of orthogonality in its natural minimal nontrivial
context of 2-dimensionality, then we can take the above property (4.2) as a {\it definition of
orthogonality, without the use of inner product}, instead of referring to the more general
condition (3.1) in Definition 1. \\

{\bf Remark 3.} \\

It is useful to note the following. For every two vectors $( a, b ) \in {\cal I}nd\,^2\,_E$,
there exists an inner product $< ~,~ >_{a,\,b}$ on $E$, such that $< a, b >_{a,\,b} \,=\, 0$,
that is, $a$ and $b$ are orthogonal with respect to $< ~,~ >_{a,\,b}$\,, or equivalently,
$( a, b ) \in {\cal O}rt\,^2\,_E ( < ~,~ >_{a,\,b} )$. This, however, does not contradict
(4.2), since in such a case we still have, see (4.1) \\

(4.3) $~~~~~~\begin{array}{l}
     \forall ~~~~ x \in E ~: \\ \\
     ~~~~~~~ x ~=~ \lambda ( a, b, x )\, a + \mu ( a, b, x )\, b
           \end{array}$ \\

with, see (2.3), (4.2) \\

(4.4) $~~~~~~ \begin{array}{l}
                 \lambda ( a, b, x ) ~=~ < a, x >_{a,\,b} / < a, a >_{a,\,b} \\ \\
                 \mu ( a, b, x ) ~=~ < b, x >_{a,\,b} / < b, b >_{a,\,b}
              \end{array} $ \\ \\

{\bf 5. One possible meaning of "Does not depend on"} \\

Here we shall specify within a rather general context
 one possible meaning of the above property "does not depend on", which was used in
(3.1) and (4.2). \\

Given a function $f : \Delta \subseteq X \times Y \longrightarrow Z$, where $X, Y$ and $Z$ are
arbitrary sets, and further given $\Gamma \subseteq \Delta$, we say that, on $\Gamma$, the
function $f$ {\it does not depend on} $x \in X$, if and only if there exists a function $g :
pr_Y ( \Gamma ) \subseteq Y \longrightarrow Z$, such that \\

(5.1) $~~~~~~ f|\,_\Gamma ~=~ g \circ pr_Y|\,_\Gamma $ \\

or equivalently, we have the commutative diagram \\

\begin{math}
\setlength{\unitlength}{0.2cm}
\thicklines
\begin{picture}(60,20)

\put(11,16){$\Delta$}
\put(28,18){$f$}
\put(16,16.5){\vector(1,0){28.5}}
\put(47,16){$~Z$}
\put(0,9.5){$(5.2)$}
\put(12,5.5){\vector(0,1){9}}
\put(9,9.5){$\subseteq$}
\put(48,5.5){\vector(0,1){9}}
\put(50,9.5){$g$}
\put(10,2){$~~ \Gamma$}
\put(16,2.6){\vector(1,0){28.5}}
\put(46.5,2){$pr_Y ( \Gamma )$}
\put(27.5,-0.5){$pr_Y|\,_\Gamma$}
\put(14.5,4.5){\vector(3,1){30}}
\put(26.5,11){$~f|\,_\Gamma$}

\end{picture}
\end{math} \\

where $pr_Y : X \times Y \ni ( x, y ) \longmapsto y \in Y$ is the usual projection mapping. \\

The above concept can easily be adapted to the case of $m \geq 2$ functions \\

(5.3) $~~~~~~ f_1,~.~.~.~, f_m : \Delta \subseteq X_1 \times ~.~.~.~ \times X_m \times Y
                                                                   ~\longrightarrow~ Z $ \\

where $X_1, ~.~.~.~ , X_m, Y$ and $Z$ are arbitrary sets. Indeed, given $\Gamma \subseteq
\Delta$, we say that, on $\Gamma$, each function $f_i$, with $1 \leq i \leq m$, {\it does not
depend on} $x_j \in X_j$, with $1 \leq j \leq m,~ j \neq i$, if and only if there exist
functions \\

(5.4) $~~~~~~ g_i : pr_{X_i \times Y} ( \Gamma ) \subseteq X_i \times Y
                                    ~\longrightarrow~ Z,~~~ 1 \leq i \leq m $ \\

such that \\

(5.5) $~~~~~~ f_i|\,_\Gamma ~=~ g_i \circ pr_{X_i \times Y}|\,_\Gamma,~~~1 \leq i \leq m $ \\

Returning to the Definition in section 2, we note that there we have \\

(5.6) $~~~~~~ \begin{array}{l}
                 X_1 = ~.~.~.~ = X_m = Y = V,~~ Z = \mathbb{R} \\ \\
                 \Delta ~=~ \bigcup_{ ( a_1, ~.~.~.~ , a_m ) \in~ {\cal I}nd\,^m\,_V}~
                   \{~ ( a_1, ~.~.~.~ , a_m ) ~\} \times span ( a_1, ~.~.~.~ , a_m ) \\ \\
                 \Gamma ~=~ \bigcup_{ ( a_1, ~.~.~.~ , a_m ) \in~ {\cal O}rt\,^m\,_V}~
                   \{~ ( a_1, ~.~.~.~ , a_m ) ~\} \times span ( a_1, ~.~.~.~ , a_m ) \\ \\
                 f_i ( a_1, ~.~.~.~ , a_m, x ) ~=~ \lambda_i ( a_1, ~.~.~.~ , a_m, x ), \\
                 ~~~~~~~~~~~~~~~~~~~~~~~~~~~~~~~~\mbox{for}~~ 1 \leq i \leq m,~
                                        ( a_1, ~.~.~.~ , a_m, x ) \in \Delta \\ \\
                  g_i ( a_i, x ) ~=~ < a_i, x > / < a_i, a_i >, \\
                 ~~~~~~~~~~~~~~~~~~~~~~~~~~~~~~~~~\mbox{for}~~ 1 \leq i \leq m,~
                                          a_i \in X_i,~a_i \neq 0,~ x \in Y
              \end{array} $ \\

It follows that \\

(5.7) $~~~~~~ {\cal O}rt\,^m\,_V ~=~ pr_{X_1 \times ~.~.~. \times X_m}\, \Gamma $ \\

And in the particular case when \\

(5.8) $~~~~~~ dim~ V ~=~ m $ \\

we obtain the simpler forms \\

(5.9) $~~~~~~ \begin{array}{l}
                 \Delta ~=~ {\cal I}nd\,^m\,_V~ \times V \\ \\
                 \Gamma ~=~ {\cal O}rt\,^m\,_V~ \times V \\ \\
                 {\cal O}rt\,^m\,_V ~=~ pr_{V^m}\, \Gamma
              \end{array} $ \\

Furthermore, through any vector space isomorphism between $V$ and $\mathbb{R}^m$, we obtain \\

(5.10) $~~~~~~ {\cal I}nd\,^m\,_V ~~~\mbox{is open in}~~ V^m,~~~~
                               \Delta ~~~\mbox{is open in}~~ V^{m + 1} $ \\

Obviously, for $m = 2$, the relations (5.9) and (5.10) apply as well to the situation in (4.1)
and (4.2). \\ \\

{\bf 6. Maximality} \\

With the notation in (5.3) - (5.5), let $\Gamma \subseteq \Delta$, and let us suppose that, on
$\Gamma$, each function $f_i$, with $1 \leq i \leq m$, does not depend on $x_j \in X_j$, with
$1 \leq j \leq m,~ j \neq i$. Then obviously, the same holds on every subset $\Gamma^\prime
\subseteq \Gamma$. \\

Let now $\Gamma_\alpha \subseteq \Delta$, with $\alpha \in A$, be a family of subsets which is
totally ordered by inclusion. Further, let us suppose that on each $\Gamma_\alpha$, each of
the functions $f_i$, with $1 \leq i \leq m$, does not depend on $x_j \in X_j$, with $1 \leq
j \leq m,~ j \neq i$. The according to (5.4), (5.5), there exist functions \\

(6.1) $~~~~~~ g_{\alpha, \,i} : pr_{X_i \times Y} ( \Gamma_\alpha ) \subseteq X_i \times Y
                                    ~\longrightarrow~ Z,~~~ \alpha \in A,~~1 \leq i \leq m $ \\

such that \\

(6.2) $~~~~~~ f_i|\,_{\Gamma_\alpha} ~=~
            g_{\alpha, \,i} \circ pr_{X_i \times Y}|\,_{\Gamma_\alpha},~~~
                                              \alpha \in A, ~~1 \leq i \leq m $ \\

If we consider now the subset \\

(6.3) $~~~~~~ \Gamma ~=~ \bigcup_{\alpha \in A}~ \Gamma_\alpha ~\subseteq~ \Delta $ \\

then for every $1 \leq i \leq m$, we obviously have \\

(6.4) $~~~~~~ pr_{X_i \times Y} ( \Gamma ) ~=~
       \bigcup_{\alpha \in A}~ pr_{X_i \times Y}( \Gamma_\alpha ) $ \\

hence in view of (6.1), (6.2), there exists a function \\

(6.5) $~~~~~~ g_i : pr_{X_i \times Y} ( \Gamma ) \subseteq X_i \times Y
                                    ~\longrightarrow~ Z $ \\

such that \\

(6.6) $~~~~~~ f_i|\,_\Gamma ~=~
            g_i \circ pr_{X_i \times Y}|\,_\Gamma $ \\

In view of the above and the Zorn lemma, we obtain \\

{\bf Lemma} \\

Given \\

(6.7) $~~~~~~ f_1,~.~.~.~, f_m : \Delta \subseteq X_1 \times ~.~.~.~ \times X_m \times Y
                                                                   ~\longrightarrow~ Z $ \\

and a nonvoid subset $\Gamma \subseteq \Delta$, such that, on $\Gamma$, each function $f_i$,
with $1 \leq i \leq m$, does not depend on $x_j \in X_j$, with $1 \leq j \leq m,~ j
\neq i$. \\
Then there exists \\

(6.8) $~~~~~~ \mbox{a maximal subset}~~~~ \bar \Gamma \subseteq \Delta,~~~
                                 \mbox{with}~~~ \bar \Gamma \supseteq \Gamma $ \\

such that, on $\bar \Gamma$, each function $f_i$, with $1 \leq i \leq m$, does not depend on
$x_j \in X_j$, with $1 \leq j \leq m,~ j \neq i$. \\

{\bf Remark 4.} \\

A similar Lemma obviously holds for the particular situation in (5.1), (5.2). \\

{\bf Theorem} \\

$\Gamma$ as defined in (5.6) is maximal. \\

{\bf Proof} \\

It follows easily. \\ \\

{\bf 7. The self-referentiality or circularity of orthogonality} \\

In section 5 one possible precise meaning of "does not depend on", a concept used in (3.1) and (4.2),
was presented. However, in view of the expressions of $\Gamma$ and ${\cal O}rt\,^m\,_V$ in (5.6),
(5.7), or for that matter, in (5.9), it is obvious that a self-referentiality or circularity
appears with respect to the concept of orthogonality when one defines it - as in (3.1), (3.2)
or (4.2) - {\it without} any inner product. \\

In this regard, and in view of recent developments in what is called {\it non-wellfounded set
theory}, see Barwise \& Moss, Forti \& Honsell, Acz\'{e}l, and the literature cited there, such
self-referential or circulatory definitions are acceptable. \\

A detailed application of non-welfounded set theory to the above problem of orthogonality will be presented in a number of
subsequent papers. \\ \\

{\bf Reference} \\

Acz\'{e}l P : Non-Well Founded Sets. Lecture Notes No. 14, CSLI Publications, 1988 \\

Amir, D : Characterizations of Inner Product Spaces. Birkh\"{a}user, Basel, 1986 \\

Barwise J, Moss L : Vicious Circles, on the Mathematics of Non-Wellfounded Phenomena.
Lecture Notes No. 60, CSLI Publications,1996 \\

Forti M, Honsell F : Set Theory with Free Construction Principles. Annali Scuola Normale
Superiore de Pisa, Classe di Scienze, 1983, No. 10, 493-522

\end{document}